\documentclass[12pt,fleqn]{article}
\usepackage{graphicx}

\usepackage{latexsym}

\usepackage{amsmath}
\usepackage{amsthm}
\usepackage{amssymb}
\usepackage{amsfonts}

\numberwithin{equation}{section}

\begin{document}

\newcommand{\rf}[1]{(\ref{#1})}
\newcommand{\rff}[2]{(\ref{#1}\ref{#2})}

\newcommand{\ba}{\begin{array}}
\newcommand{\ea}{\end{array}}

\newcommand{\be}{\begin{equation}}
\newcommand{\ee}{\end{equation}}

\newcommand{\const}{{\rm const}}
\newcommand{\ep}{\varepsilon}
\newcommand{\Cl}{{\cal C}}
\newcommand{\rr}{{\vec r}}
\newcommand{\ph}{\varphi}
\newcommand{\R}{{\mathbb R}}  
\newcommand{\C}{{\mathbb C}}  
\newcommand{\T}{{\mathbb T}}  
\newcommand{\Z}{{\mathbb Z}}  

\newcommand{\e}{{\bf e}}

\newcommand{\m}{\left( \ba{c}}
\newcommand{\ema}{\ea \right)}
\newcommand{\mm}{\left( \ba{cc}}
\newcommand{\miv}{\left( \ba{cccc}}

\newcommand{\av}[1]{\mbox{$\langle #1 \rangle$}} 
\newcommand{\scal}[2]{\mbox{$\langle #1 \! \mid #2 \rangle $}}
\newcommand{\ods}{\par \vspace{0.5cm} \par}
\newcommand{\dis}{\displaystyle }
\newcommand{\mc}{\multicolumn}
\newcommand{\id}{{\rm id}}

\newtheorem{prop}{Proposition}[section]
\newtheorem{Th}[prop]{Theorem} 
\newtheorem{lem}[prop]{Lemma}
\newtheorem{rem}[prop]{Remark}
\newtheorem{cor}[prop]{Corollary}
\newtheorem{Def}[prop]{Definition}
\newtheorem{open}{Open problem}
\newtheorem{ex}{Example}
\newtheorem{exer}{Exercise}

\newenvironment{Proof}{\par \vspace{2ex} \par
\noindent \small {\it Proof:}}{\hfill $\Box$ 
\vspace{2ex} \par }

\title{\bf 
New definitions of exponential, hyperbolic and trigonometric functions on time scales}
\author{
 {\bf Jan L.\ Cie\'sli\'nski}\thanks{\footnotesize
 e-mail: \tt janek\,@\,alpha.uwb.edu.pl}
\\ {\footnotesize Uniwersytet w Bia{\l}ymstoku,
Wydzia{\l} Fizyki}
\\ {\footnotesize ul.\ Lipowa 41, 15-424
Bia{\l}ystok, Poland}
}

\date{}

\maketitle

\begin{abstract}
We propose two new definitions of the exponential function on time scales. The first definition is based on the Cayley transformation while the second one is a natural extension of exact discretizations. Our eponential functions map the imaginary axis into the unit circle. Therefore, it is possible to define hyperbolic and trigonometric functions on time scales  in a standard way. The resulting functions preserve most of the qualitative properties of the corresponding continuous functions. In particular, Pythagorean trigonometric identities hold exactly on any time scale. Dynamic equations satisfied by Cayley-motivated functions have a natural similarity to the corresponding diferential equations. 
The exact discretization is less convenient as far as dynamic equations and differentiation is concerned. 
\end{abstract}

\ods

{\it MSC 2010:} 33B10; 26E70; 34N05; 65L12

{\it Key words and phrases:} time scales, measure chains, exponential function, trigonometric functions, hyperbolic functions, Cayley transform,  first and second order dynamic equations, exact discretization

\pagebreak

\section{Introduction}

The main goal of the time scales approach is a unification of the differential and difference calculus \cite{Hi,Hi2,ABOP}. 
There are many papers and books where notions and notation concerning time scales are explained in detail \cite{Hi,ABOP,BP-I}, see also \cite{BG-Lap}. 
A time scale $\T$ is defined as an arbitrary closed subset of $\R$ \cite{Hi,ABOP}. The forward jump operator $\sigma$ is defined as $\sigma (t) = \inf \{ s \in \T : s > t \}$ (and we assume $\inf \emptyset = \sup \T$). We usually denote $t^\sigma := \sigma (t)$. 

A point $t \in \T$ is called right-dense iff $t^\sigma  = t$ and right-scattered iff $t^\sigma  > t$. Similarly one can define the backward jump-operator $\rho$, left-dense points and left-scattered points \cite{Hi}. Sometimes we have to exclude from considerations the left-scattered maximum of $\T$ (if it exists). $\T$ minus the left-scattered maximum will be denoted by $\T^\kappa$ (if such maximum does not exist then $\T^\kappa = \T$). 

{\it Graininess} \ $\mu = \mu (t)$ is defined, for $t \in \T^\kappa$, as $\mu (t) = t^\sigma - t$. 
 The delta derivative (an  analogue of the derivative for functions defined on time scales) is defined by
\be  \label{delta}
   x^\Delta (t) := \lim_{\stackrel{\dis s\rightarrow t}{s \neq \sigma (t)}} \frac{ x (\sigma (t) ) - x (s) }{\sigma (t) - s } \ , \qquad (t \in \T^\kappa)  \ .
\ee
If $t$ is right-dense, then $x^\Delta (t) = \dot x (t)$. If $t$ is right-scattered, then  $x^\Delta$ is the corresponding  difference quotient. 

An important notion is rd-continuity. A function is said to be righ-dense continuous (rd-continuous), if it is continuous at right-dense points and at left-dense points there exist finite left-hand limits. 
The graininess $\mu$ is rd-continuous but, in general, is not continuous. 
{\it Dynamic equations} \ are time scales analogues of differential equations (i.e., they may contain delta derivatives, jump operators and sometimes the graininess of the considered time scale) \cite{ABOP,BP-I}. 

In this paper we consider the problem of defining elementary (or special) functions on time scales. We propose a new approach to  exponential, hyperbolic and trigonometric functions. 
We point out that the definition of such functions is not unique. 
We suggest new definitions (improved exponential, hyperbolic and trigonometric functions) based on the Cayley transformation. The new functions preserve  more properties of their continuous counterparts in comparison to the existing definitions. In particular, our exponential function maps the imaginary axis into the unit circle and trigonometric functions satisfy the Pythagorean identity. Dynamic equations satisfied by our improved functions have a natural similarity to the corresponding differential equations.  

We also propose the notion of {\it exact} time scales analogues of special functions and we identify Hilger's definition of trigonometric functions \cite{Hi-spec} with exact analogues of these functions.  Exact discretizations of differential equations \cite{Po,Mic,Ag} are intended as a way to connect smooth and discrete case, but in an apparently different way than the time scales calculus. In this paper we incorporate the exact discretizations in the framework of the time scales approach.
We discuss exact exponential, hyperbolic and trigonometric functions on constant time scales. An extension of these results on more general time scales seems to be rather difficult.

\section{Survey of existing definitions}

The exponential function on time scales has been introduced by Hilger \cite{Hi}, and his definition seems to be commonly accepted \cite{ABOP,BP-I}. 
A different situation has place in the case of hyperbolic and trigonometric functions, where 
Hilger's approach \cite{Hi-spec} differs from Bohner-Peterson's approach \cite{BP12}. 

\subsection{Exponential function}

Hilger defined the exponential function as follows:
\be  \label{exp-Hi}
   e_\alpha (t, \tau) := \exp \left( \int_\tau^t \xi_{\mu (s)} ( \alpha (s)) \Delta s \right) \ ,  \qquad 
e_\alpha (t) := e_\alpha (t, 0) \ , 
\ee
where 
\be
  \xi_h (z) :=  \frac{1}{h} \log (1 + z h) \quad ({\rm for} \ \  h > 0)  \quad {\rm and} \ \ \xi_0 (z) := z \ . 
\ee
This definition applies to the so called $\mu$-{\it regressive} functions $\alpha = \alpha (t)$, i.e., those satisfying
\be
  1 + \mu (t) \alpha (t) \neq 0 \quad {\rm for \ all} \quad t\in \T^\kappa  \ .
\ee
Such functions are usually called regressive, but we reserve this name for another class of functions, see Definition~\ref{Def-regr}. 

In the constant discrete case ($\T = \ep \Z$, $\alpha = \const$) we have 
\be
  e_\alpha (t) = (1 + \alpha \ep)^{\frac{t}{\ep}} \ ,
\ee
and in the case $\T = \R$ we have
\be
 e_\alpha (t) = \exp \int_0^t \alpha (\tau) d \tau \ .
\ee 

\begin{Th}[\cite{Hi,BP12}]
If $\alpha, \beta : \T \rightarrow \C$ are $\mu$-regressive and rd-continuous, then 
the following properties hold: 
\begin{enumerate}
\item $e_\alpha (t^\sigma, t_0) = (1 + \mu (t) \alpha (t) ) \ e_\alpha (t, t_0)$ \ , 
\item $( e_\alpha (t, t_0) )^{-1} = e_{\ominus^\mu \alpha} (t, t_0) $ \ ,
\item $e_\alpha (t, t_0) \ e_\alpha (t_0, t_1) = e_\alpha (t, t_1)$ \ , 
\item $e_\alpha (t, t_0) \ e_\beta  (t, t_0) =  e_{\alpha\oplus^\mu \beta} (t, t_0) $ \ ,
\end{enumerate}
where \ $\alpha \oplus^\mu \beta := \alpha + \beta + \mu \alpha \beta$ \ and \ $\ominus^\mu \alpha := \frac{- \alpha}{1 + \mu \alpha}$. 
\end{Th}

The addition $\oplus^\mu$ is usually denoted by $\oplus$. However, we reserve the notation $\oplus$ for another addition, see Definition~\ref{def-oplus}. The exponential function \rf{exp-Hi} solves the Cauchy problem: $x^\Delta = \alpha x$, $x (0) = 1$, see \cite{Hi}. 

Similar considerations (with the delta derivative replaced by the nabla derivative) lead to the 
nabla exponential function \cite{ABEPT}. In the case $\T = \ep \Z$ the nabla exponential function is given by
\be
 {\hat e}_\alpha (t) = (1 - \alpha \ep)^{- \frac{t}{\ep}} 
\ee
and solves the Cauchy problem: $x^\nabla = \alpha x$, $x (0) = 1$. The existence of the above two definitions reflects the duality between delta and nabla calculus. 
A linear combination of delta nad nabla derivatives (the diamond-$\alpha$ derivative \cite{SFHD}) leads to another definition of the exponential function (the so-called diamond-alpha exponential function \cite{MT}).

\subsection{Hilger's approach to hyperbolic and trigonometric functions}

The first approach to hyperbolic functions has been proposed by Hilger \cite{Hi-spec},  
\be  \label{hyp-Hi}
   \cosh_\alpha (t) = \frac{e_\alpha (t) + e_{\ominus^\mu\alpha} (t) }{2} \ ,  \qquad 
\sinh_\alpha (t) = \frac{e_\alpha (t) - e_{\ominus^\mu\alpha} (t) }{2}  \ ,
\ee
where $\alpha$ is $\mu$-regressive. Among its advantages we have the identity   
\be  \label{Pyt-hyp}
\cosh_\alpha^2 (t) - \sinh_\alpha^2 (t) = 1 \ . 
\ee
Delta derivatives of these functions are linear combinations of both hyperbolic functions, e.g., 
\be
 \cosh_\alpha^\Delta (t) = \frac{\alpha + [\ominus^\mu \alpha]}{2} \cosh_\alpha (t) + \frac{\alpha - [\ominus^\mu \alpha]}{2} \sinh_\alpha (t) \ .
\ee
In the constant discrete case ($\T = \ep \Z$, $\alpha = \const$) we have 
\be \ba{l} \dis
   \cosh_\alpha (t) =  
\frac{ (1 + \alpha \ep)^{\frac{t}{\ep}} + (1 + \alpha \ep)^{-\frac{t}{\ep}} }{2} \ ,  \\[3ex] \dis
\sinh_\alpha (t) = 
\frac{ (1 + \alpha \ep)^{\frac{t}{\ep}} - (1 + \alpha \ep)^{-\frac{t}{\ep}} }{2} \ .
\ea \ee
The hyperbolic functions \rf{hyp-Hi} evaluated at imaginary $\alpha$ are not real. Thus the definition \rf{hyp-Hi} cannot be extended on trigonometric functions by substituting $i \omega$ for $\alpha$. 
In order to treat trigonometric functions Hilger introduced  another map $\alpha \rightarrow  \overset{\circ}{\iota} \omega$, see \cite{Hi-spec}. In the constant case (i.e., $\mu (t) = \const$ and $\omega = \const$) the final result is very simple
\be
\cos_\omega (t) = \cos \omega t \ , \quad \sin_\omega (t) = \sin\omega t \ .
\ee
In fact, this is a particular case of the {\it exact discretization} (see Section~\ref{sec-exact}), although Hilger's motivation seems to be different. Exact discretizations have many advantages (compare \cite{Ci-oscyl})  but, unfortunatelly, their delta derivatives are quite complicated. For instance, at right scattered points we have:
\be  \label{delta-sin-ex}
 (\sin\omega t)^\Delta = \frac{\sin\omega\mu}{\mu} \cos\omega t + \frac{\cos\omega\mu - 1}{\mu} \sin\omega t \ ,   
\ee 
compare \cite{Hi-spec}, see also Section~\ref{sec-exact}.

\subsection{Bohner-Peterson's approach to hyperbolic and tri\-gonometric functions}

The second approach has been proposed by Bohner and Peterson \cite{BP12,BP-exp}
\be  \label{hyp-BP}
\ba{l} \dis 
   \cosh_\alpha (t) = \frac{e_\alpha (t) + e_{-\alpha} (t) }{2} \ ,  
\\[2ex] \dis
\sinh_\alpha (t) = \frac{e_\alpha (t) - e_{-\alpha} (t) }{2}  \ ,
\ea \ee
where $\alpha$ is $\mu$-regressive. The hyperbolic functions defined by \rf{hyp-BP} satisfy 
\be
 \cosh_\alpha^\Delta (t) = \alpha \sinh_\alpha (t) \ , \qquad 
\sinh_\alpha^\Delta (t) = \alpha \cosh_\alpha (t) \ , 
\ee
The identity \rf{Pyt-hyp} is not valid. Instead we have
\be  \label{Pyt-BP-hyp}
\cosh_\alpha^2 (t) - \sinh_\alpha^2 (t) = e_{-\mu\alpha^2} (t) \ .
\ee
Bohner and Peterson define trigonometric functions in a natural way, evaluating hyperbolic functions at the imaginary axis
\be
\cos_\omega (t) = \cosh_{i \omega} (t) \ , \qquad  i \sin_\omega (t) = \sinh_{i \omega} (t) \ .
\ee
Then,
\be  \label{Pyt-BP-trig} 
\cos_\omega^2 (t) + \sin_\omega^2 (t) = e_{\mu\omega^2} (t) \ .
\ee
In the constant discrete case ($\T = \ep \Z$, $\alpha = \const$, $\omega = \const$) we have
\be  
\ba{l} \dis 
   \cosh_\alpha (t) = 
\frac{ (1 + \alpha \ep)^{\frac{t}{\ep}} + (1 - \alpha \ep)^{\frac{t}{\ep}} }{2} \ ,  \\[2ex] \dis
\sinh_\alpha (t) = 
\frac{ (1 + \alpha \ep)^{\frac{t}{\ep}} - (1 - \alpha \ep)^{\frac{t}{\ep}} }{2} \ . 
\ea \ee 
Moreover, 
\be
e_{-\ep \alpha^2} (t) = (1 - \ep^2 \alpha^2)^{\frac{t}{\ep}} \ , \quad 
e_{\ep \omega^2} (t) = (1 + \ep^2 \omega^2)^{\frac{t}{\ep}} \ .  
\ee
Therefore, 
\be
\lim_{t \rightarrow - \infty} e_{-\ep \alpha^2} (t) = \infty \ , \qquad 
\lim_{t \rightarrow \infty} e_{-\ep \alpha^2} (t) =  0 \ , 
\ee
provided that  $|\alpha \ep| < 1$. Similarly,   
\be
\lim_{t \rightarrow - \infty} e_{\ep \omega^2} (t) = 0  \ , \qquad 
\lim_{t \rightarrow \infty} e_{\ep \omega^2} (t) =  \infty \ , 
\ee
Therefore, the definition \rf{hyp-BP} leads to Pythagorean-like identities \rf{Pyt-BP-hyp}, \rf{Pyt-BP-trig} which have essentially different behaviour in the discrete and continuous case.

\section{Approach motivated by the Cayley transformation}
\label{sec-main}

In this section we present new definitions of exponential, hyperbolic and trigonometric functions and their properties. We will tentatively refer to them as `improved' functions, because they simulate the behaviour of their continuous counterparts better than the previous definitions.  Our new definitions are based on the classical Cayley transformation:
\be  \label{cay}
  z \rightarrow {\rm cay} (z,a) = \frac{1 + a z}{1- a z} \ ,
\ee
see, for instance, \cite{Is-Cay}.

\subsection{New definition of the exponential function}

In order to formulate our definition we need to redefine a notion of regressivity.  

\begin{Def}  \label{Def-regr}
The function $\alpha : \T \in \C$ is regressive if \ $\mu (t) \alpha (t) \neq \pm 2$ \ for any $t \in \T^\kappa$. 
\end{Def}

\begin{Def} 
The improved exponential function (or the Cayley-exponen\-tial function) on a time scale  is defined by
\be  \label{E-for}
  E_\alpha (t, t_0) := \exp \left( \int_{t_0}^t \zeta_{\mu (s)} ( \alpha (s)) \Delta s \right) \ , \qquad E_\alpha (t) := E_\alpha (t, 0) \ , 
\ee
where $\alpha = \alpha (t)$ is a given rd-continuous regressive function and
\be \label{zeta}
  \zeta_h (z) :=  \frac{1}{h} \log \frac{1 + \frac{1}{2} z h}{1 - \frac{1}{2} z h} \quad ({\rm for} \ \  h > 0)  \quad {\rm and} \ \ \zeta_0 (z) := z \ .
\ee
\end{Def}

Here and in what follows the logarithm is understood as a principal branch of the complex logaritm with image $[-i\pi, i\pi]$.  

\begin{lem}
If $\alpha$ is rd-continuous and regressive, then the delta-integral in \rf{E-for} exists.
\end{lem}

\begin{Proof}
The assumption of regressivity of $\alpha$ implies that the logarithm in \rf{E-for} exists (is finite) for any $t \in \T^\kappa $. Thus 
the function $t \rightarrow \zeta_\mu (t) (\alpha (t))$ has no singularities. To complete the proof we will show that   
$\zeta_\mu \circ \alpha $ is rd-continuous (which implies that it has  an antiderivative, see \cite{Hi}). At right-dense $t_0$ we have 
\be
\lim_{t\rightarrow t_0} \alpha (t) = \alpha (t_0) \ , \qquad 
\lim_{t\rightarrow t_0} \mu (t) = \mu (t_0) = 0 \ , 
\ee
because $\alpha$ and $\mu$ are continuous at right-dense points. Therefore
\be
\lim_{t\rightarrow t_0} \zeta_{\mu (t)} (\alpha (t)) = \lim_{t\rightarrow t_0} \frac{1}{\mu (t)} \log \frac{1 + \frac{1}{2} \alpha (t) \mu (t) }{1 - \frac{1}{2} \alpha (t) \mu (t) } =  \alpha (t_0) \ .
\ee
On the other hand, $\zeta_{\mu (t_0)} (\alpha (t_0)) = \zeta_0 (\alpha (t_0)) = \alpha (t_0)$. Therefore, $\zeta_\mu \circ \alpha$ is continuous at right-dense points.  At left-dense $s_0$ we denote  
\be
\alpha (s_0^-) := \lim_{t\rightarrow s_0^-} \alpha (t) \ , \qquad 
\mu (s_0^-) := \lim_{t\rightarrow s_0^-} \mu (t) = 0 \ .
\ee
In general $\alpha (s_0^-) \neq \alpha (s_0)$, $\mu (s_0^-) \neq \mu (s_0)$ but rd-continuity guarantees that all these values are finite. Then, 
\be
\lim_{t\rightarrow s_0^-} \zeta_{\mu (t)} (\alpha (t)) = \lim_{t\rightarrow s_0^-} \frac{1}{\mu (t)} \log \frac{1 + \frac{1}{2} \alpha (t) \mu (t) }{1 - \frac{1}{2} \alpha (t) \mu (t) } =  \alpha (s_0^-) \ , 
\ee
and the existence of this finite limit means that $\zeta_\mu \circ \alpha$ is rd-continuous. 
\end{Proof}

\ods

In the constant discrete case ($\T = \ep \Z$, $\alpha = \const$) we have 
\be  \label{exp-dis}
  E_\alpha (t) = \left( \frac{1 + \frac{1}{2} \alpha \ep}{1 - \frac{1}{2} \alpha \ep} \right)^{\frac{t}{\ep}} \ ,
\ee
and in the case $\T = \R$ we have, as usual, 
\be
 E_\alpha (t) = \exp \int_0^t \alpha (\tau) d \tau \ .
\ee 
The formula \rf{exp-dis} appeared earlier in different contexts, see for instance \cite{Is-Cay,Mer}. 

The new definition \rf{E-for} of the exponential function can be related to Hilger's  definition \rf{exp-Hi} with an other exponent. Namely,  
we are going to prove that $E_\alpha (t, t_0) = e_\beta (t, t_0)$   provided that
\be  \label{betal}
\beta (t) = \frac{ \alpha (t)}{1 - \frac{1}{2}{\mu (t) \alpha (t)}} \ .
\ee

\begin{Th} \label{Ee}
For any regressive, rd-continuous $\alpha = \alpha (t)$ there corresponds a unique  $\mu$-regressive, rd-continuous $\beta = \beta (t)$ (given by \rf{betal}) such that 
$E_\alpha (t, t_0) = e_\beta (t, t_0)$. 
For $\mu$-regressive, rd-continuous $\beta$ satisfying  $\mu \beta \neq -2$ there exists a unique regressive, rd-continuous $\alpha$ given by
\be  \label{albet}
   \alpha (t) = \frac{\beta (t)}{1 + \frac{1}{2} \mu (t) \beta (t)} \ ,
\ee
such that $E_\alpha (t, t_0) = e_\beta (t, t_0)$. 
\end{Th}

\begin{Proof}
$E_\alpha (t, t_0) = e_\beta (t, t_0) $ if and only if the integrands in \rf{exp-Hi} and \rf{E-for} coincide, i.e., $\xi_\mu ( \beta) = \zeta_\mu (\alpha) $. Thus $\beta (t) = \alpha (t)$ at right-dense points, and for $\mu \neq 0$ (i.e., at right-scattered points) we get
\be  
1 + \mu \beta = \frac{1 + \frac{1}{2} \mu \alpha}{1 - \frac{1}{2} \mu \alpha} \ .
\ee
Both cases lead to a single condition \rf{betal} (or, equivalently, to \rf{albet}). To complete the proof we verify that $\mu \beta = -1$ iff $\mu \alpha = - 2$. Therefore for any regressive $\alpha$ we have $\mu \beta \neq - 1$. Then, $\mu \alpha = 2$ corresponds to $\mu \beta = \pm \infty$. Thus for any $\mu$-regressive $\beta$ we have $|\mu \alpha| \neq 2$. Finally, we observe that $\mu \beta= -2$ corresponds to $\mu \alpha = \pm \infty$, and for any other values of $\beta$ the value of  $\alpha$ is uniquely determined by \rf{albet}.

To complete the proof we have to show that $\alpha$ is rd-continuous iff $\beta$ is rd-continuous. 
 At right-dense $t_0$ functions $\alpha$, $\mu$ are continuous and also  $\mu (t_0) = 0$. Therefore,  from \rf{albet} we get
\be
\lim_{t\rightarrow t_0} \alpha (t) = \lim_{t\rightarrow t_0} \beta (t)  \ , \qquad \alpha (t_0) = \beta (t_0)   \ .
\ee
Hence, $\alpha$ is continuous at $t_0$ if and only if $\beta$ is continuous at $t_0$. If $s_0$ is left-dense, then \ $\mu (t) \rightarrow 0$ \ as \ $t \rightarrow s_0^-$. As  a consequence, we have 
\be
\alpha (s_0^-) \equiv \lim_{t\rightarrow s_0^-} \alpha (t) =  \lim_{t\rightarrow s_0^-} \beta (t) \equiv  \beta (s_0^-) \ . 
\ee
Therefore,  $\alpha (s_0^-)$ exists if and only if $\beta (s_0^-)$ exists, which ends the proof. 
\end{Proof}

\subsection{Properties of the Cayley-exponential function} 

First of all, we observe a close relation between $\zeta_h$ and the Cayley transformation \rf{cay}, 
\be
  e^{h \zeta_h (z)} = {\rm cay} (z, \frac{1}{2} h) \ .
\ee
The transformation inverse to $\zeta_h$ is given by
\be  \label{zth}
 z \equiv \zeta_h^{-1} (\zeta) = \frac{2}{h} \tanh \frac{h \zeta}{2} \quad (h \neq 0) \ , \qquad \zeta_0^{-1} (\zeta) = \zeta \ .
\ee
In particular, 
\be
z =  \frac{2 i}{h} \tan \frac{ h \omega}{2} \qquad {\rm for} \quad   
\zeta = i \omega \ .
\ee
Therefore, $\zeta_h^{-1}$ maps the open segment $(-\frac{\pi i}{h}, \frac{\pi i}{h}) \subset i \R$ onto  $i \R$, and $\R$ is mapped onto the real segment $(-\frac{2}{h}, \frac{2}{h})$.

\begin{cor}
$\zeta_h$ maps the imaginary axis onto the segment $(-\frac{\pi i}{h}, \frac{\pi i}{h}) \subset i \R$. 
\end{cor}

\ods
\begin{lem} Denoting $\zeta = \gamma + i \eta$ and taking into account \rf{zth}, we have
\be \ba{l} \dis
|z | < \frac{2}{h}  \quad \Longleftrightarrow \quad \cos \eta h > 0 \ , \\[2ex]\dis
|z | = \frac{2}{h} \quad \Longleftrightarrow \quad \cos \eta h = 0 \ , \\[2ex]\dis
|z | > \frac{2}{h} \quad \Longleftrightarrow \quad \cos \eta h < 0 \ .
\ea \ee
\end{lem}

\begin{Proof} We compute:
\be
  \tanh \frac{h \zeta}{2} = \frac{ e^{h \gamma + i h \eta} - 1  }{ e^{h \gamma + i h \eta} + 1  } = \frac{ e^{h \gamma} \cos h \eta - 1 + i e^{h \gamma} \sin h \eta }{ e^{h \gamma} \cos h \eta + 1 + i e^{h \gamma} \sin h \eta } \ ,
\ee
\be
 \left| \tanh \frac{h \zeta}{2} \right| = \sqrt{ \frac{ e^{2 h \gamma} - 2 e^{h \gamma} \cos h \eta + 1 }{ e^{2 h \gamma} + 2 e^{h \gamma} \cos h \eta + 1 } } \ .
\ee
To complete the proof it is enough to notice that $e^{h \gamma} > 0$ and \ $\frac{h z}{2} = \tanh\frac{h \zeta}{2} $. 
\end{Proof}

\ods
\begin{cor}
$\zeta_h$ maps the disc \ $|z| < \frac{2}{h}$ \ onto the strip \ $ - \frac{\pi}{2 h} < \eta < \frac{\pi}{2 h} $. 
\end{cor}

\ods

\begin{Def}  \label{def-oplus}
Given a time scale $\T$ and two functions $\alpha, \beta: \T \rightarrow \C$, we define
\be  \label{oplusmu}
 \alpha \oplus \beta := \frac{\alpha + \beta}{1 + \frac{1}{4} \mu^2  \alpha \beta } \ .
\ee
where $\mu = \mu (t)$ is the graininess of $\T$. 
\end{Def}

\begin{lem}  \label{lem-zeta}
The function $\zeta_\mu$ has the following properties: 
\be  
\overline{\zeta_\mu (\alpha)} = \zeta_\mu (\bar \alpha) \ , \quad \zeta_\mu (-\alpha) = - \zeta_\mu (\alpha) , \quad \zeta_{-\mu} (\alpha) = \zeta_\mu (\alpha)  \ , 
\ee
\be 
\zeta_\mu (\alpha) + \zeta_\mu (\beta) = \zeta_\mu (\alpha\oplus\beta) ,
\ee
where bar denotes the complex conjugate and we assume (in order to avoid infinities) 
$\frac{1}{2} \mu \alpha \neq - 1$, $\frac{1}{2} \mu \beta \neq -1 $ and  $\frac{1}{2} \mu \beta \neq - \left( \frac{1}{2} \mu \alpha \right)^{-1}$. 
\end{lem}

\begin{Proof} The function $\zeta_\mu = \zeta_\mu (\alpha)$ is analytic with respect to $\alpha$ (and $\mu$ is real). Hence the first property follows. Other properties can be shown by direct calculation:  
\be \ba{l} \dis  \label{proof-gamma}
 \zeta_\mu (-\alpha) = \frac{1}{\mu} \log \frac{1 - \frac{1}{2} \alpha \mu}{1 + \frac{1}{2} \alpha \mu} = - \frac{1}{\mu} \log \frac{1 + \frac{1}{2} \alpha \mu}{1 - \frac{1}{2} \alpha \mu} = - \zeta_\mu (\alpha) \ ,
\\[3ex]\dis
 \zeta_{-\mu} (\alpha) = - \frac{1}{\mu} \log \frac{1 - \frac{1}{2} \alpha \mu}{1 + \frac{1}{2} \alpha \mu} =  \frac{1}{\mu} \log \frac{1 + \frac{1}{2} \alpha \mu}{1 - \frac{1}{2} \alpha \mu} =  \zeta_\mu (\alpha) \ , 
\\[3ex]\dis
\zeta_\mu (\alpha) + \zeta_\mu (\beta) = \frac{1}{\mu} \log \frac{ 1+ \frac{1}{4} \mu^2 \alpha \beta + \frac{1}{2} \mu (\alpha +\beta) }{1+ \frac{1}{4} \mu^2 \alpha \beta - \frac{1}{2} \mu (\alpha +\beta)} = \zeta_\mu (\alpha\oplus\beta) \ ,
\ea \ee
provided that $\frac{1}{4} \mu^2 \alpha \beta \neq - 1$. 
\end{Proof}

\ods

\begin{Th} \label{Th-exp-properties}
If $\alpha, \beta : \T \rightarrow \C$ are regressive and rd-continuous, then 
the following properties hold: 
\begin{enumerate}
\item $\displaystyle E_\alpha (t^\sigma, t_0) = \frac{1 + \frac{1}{2} \mu (t) \alpha (t)}{1 - \frac{1}{2} \mu (t) \alpha (t) } \ E_\alpha (t, t_0)$ \ , 
\item $( E_\alpha (t, t_0) )^{-1} = E_{-\alpha} (t, t_0) $ \ ,
\item $ \overline{ E_\alpha (t, t_0)} = E_{\bar \alpha} (t, t_0)$ \ , 
\item $E_\alpha (t, t_0) \ E_\alpha (t_0, t_1) = E_\alpha (t, t_1)$ \ , 
\item $E_\alpha (t, t_0) \ E_\beta  (t, t_0) =  E_{\alpha\oplus \beta} (t, t_0) $ \ ,
\end{enumerate}
where we use a standard notation $t^\sigma \equiv  \sigma (t)$.

\end{Th}

\begin{Proof}
It is sufficient to prove the first property for right-scattered points ($t^\sigma > t$). 
\be \ba{l} \dis  \label{exp11}
E_\alpha (t^\sigma, t_0) = 
\exp \left( \int_{t}^{\sigma (t)} \zeta_{\mu (\tau) } (\alpha (\tau)) \Delta\tau  \right)  \exp \left( \int_{t_0}^t \zeta_{\mu (\tau) } (\alpha (\tau)) \Delta\tau  \right)
\ea \ee
Then, using \rf{zeta}, we get
\[
\int_{t}^{\sigma (t)} \zeta_{\mu (\tau) } (\alpha (\tau)) \Delta\tau  = \mu (t) \zeta_{\mu (t) } (\alpha (t)) = \log \frac{1 + \frac{1}{2} \mu (t) \alpha (t) }{ 1 - \frac{1}{2} \mu (t) \alpha (t) } \ , 
\] 
and substituting it into \rf{exp11} we get the first property.  The second property follows directly from $\zeta_\mu (-\alpha) = - \zeta_\mu (\alpha)$, see Lemma~\ref{lem-zeta}. Indeed, 
\[  \dis 
E_\alpha^{-1} (t, t_0) = 
\exp \left( - \int_{t_0}^t \zeta_{\mu (\tau) } (\alpha (\tau)) \Delta\tau  \right) = \exp  \int_{t_0}^t \zeta_{\mu (\tau) } (- \alpha (\tau)) \Delta\tau  = E_{-\alpha} (t, t_0) .
\]
The third property follows directly from analycity of the exponential function and from Lemma~\ref{lem-zeta}. We recall that $t \in \T \subset \R$.  Indeed,  
\[
\overline{ E_\alpha (t, t_0)} = 
 \exp  \int_{t_0}^t \overline{ \zeta_\mu (\tau) (\alpha (\tau))} \Delta\tau         =  \exp  \int_{t_0}^t \zeta_\mu (\tau) ( \overline{\alpha (\tau) } ) \Delta\tau      =  E_{\bar \alpha} (t, t_0) \ . 
\]
The fourth property is derived in a straightforward way:
\[
E_\alpha (t, t_0) \ E_\alpha (t_0, t_1) = \exp \left( \int_{t_0}^{t} \zeta_{\mu (t)} (\alpha (t)) \Delta t + \int_{t_1}^{t_0} \zeta_{\mu (t)} (\alpha (t)) \Delta t \right) = E_\alpha (t, t_1)  .  
\]
Finally, 
\[
E_\alpha (t, t_0) \ E_\beta (t, t_0) = \exp \int_{t_0}^{t} \left( \zeta_{\mu (t)} (\alpha (t)) + \zeta_{\mu (t)} (\beta (t)) \right) \Delta t = E_{\alpha\oplus\beta} (t, t_0)  ,  
\]
where we took into account Lemma~\ref{lem-zeta}. 
\end{Proof}
\ods

The formula \rf{oplusmu} is identical with the Lorentz velocity transformation of the special theory of relativity (the role of the speed of light $c$ is played by $\frac{2}{\mu}$). 
Denoting 
\be
\alpha':= \frac{1}{2} \mu \alpha \ , \quad 
\beta' = \frac{1}{2} \mu \beta \ .
\ee
we can rewrite the formula \rf{oplusmu} in a simpler form
\be \label{oplus'}
   \alpha'\oplus\beta' = \frac{\alpha' + \beta'}{1 + \alpha' \beta'} \ .
\ee
\ods

\begin{lem}
If $\alpha'$ and $\beta'$ are real functions on $\T$ and $\alpha'\oplus\beta'$ is given by \rf{oplus'}, then 
\be \label{abg<1} 
|  \alpha' | < 1 \ \ {\rm and} \ \ 
|  \beta' | < 1  \quad \Longrightarrow \quad 
|  \alpha'\oplus\beta' | < 1 \ , 
\ee 
\be  \label{abg=1}
|  \alpha'\oplus\beta' | = 1  \quad \Longleftrightarrow \quad 
|  \alpha' | = 1  \ \ {\rm or}  \ \ |  \beta' | = 1 \ .
\ee
\end{lem}

\begin{Proof} Using \rf{oplus'} we compute
\be
1 - (\alpha'\oplus\beta')^2 = 
\frac{(1-(\alpha')^2)(1 - (\beta')^2)}{(1 + \alpha' \beta')^2 } \ ,
\ee
which immediately yields \rf{abg=1}. Then, 
$|  \alpha' | < 1$ and $|  \beta' | < 1$ imply that the righ-hand side is  positive. Hence, we have \rf{abg<1}. 
\end{Proof}

\begin{Def}  \label{Def-posit}
The function $\alpha : \T \rightarrow \R$ is called positively regressive if for all $t \in \T^\kappa $ we have \ $ |\alpha (t) \mu (t) | < 2$. 
\end{Def}

\begin{Th}
If $\alpha : \T \rightarrow \R$ is rd-continuous and positively regressive, then the exponential function $E_\alpha$ is positive (i.e., $E_\alpha (t) > 0$ for all \ $t \in \T$). 
\end{Th}

\begin{Proof} If $\alpha$ is real and positively regressive, then (for any $t \in \T^\kappa$) we have
\be
  \frac{1 + \frac{1}{2} \mu (t) \alpha (t) }{1 - \frac{1}{2} \mu (t)  \alpha (t) } > 0 \ .
\ee
Thus $\zeta_{\mu (t)} (\alpha (t))$ is real for any $t \in \T^\kappa$ and, as a consequence, the exponential function is positive.
\end{Proof}

The regressivity condition $|\alpha (t) \mu (t) | < 2$ is automatically satisfied at  right-dense points. At right-scattered points we have
\be
   E_\alpha (t^\sigma) - E_\alpha (t) = \mu (t) \alpha (t) \frac{ E_\alpha (t^\sigma) + E_\alpha (t) }{2}  \ .
\ee
Therefore, the condition $| \alpha (t) \mu (t)  | < 2$ is equivalent to 
\be
   | E_\alpha (t^\sigma) - E^\alpha (t) | < | E_\alpha (t^\sigma) + E^\alpha (t) | 
\ee
In the real case it means that $E_\alpha (t^\sigma)$ and $E_\alpha (t)$ have the same sign.  

\ods
\begin{Th}
The set of positively regressive real functions ${\cal R}_+$  is an abelian group with respect to the addition $\oplus$.
\end{Th}

\begin{Proof} The formula \rf{oplusmu} obviously yields $\alpha\oplus\beta = \beta\oplus\alpha$. The element inverse to $\alpha$ (given simply by $\ominus\alpha = - \alpha$) always exists in ${\cal R}_+$. 
Therefore, it is sufficient to show that ${\cal R}_+$ is closed with respect to $\oplus$. 
Taking into account \rf{oplusmu} we have: 
\be
\frac{1 + \frac{1}{2} \mu \alpha \oplus \beta}{1 - \frac{1}{2} \mu \alpha \oplus \beta} =  \frac{ (1 + \frac{1}{2} \mu \alpha)(1 + \frac{1}{2} \mu  \beta)  }{ (1 - \frac{1}{2} \mu \alpha)(1 - \frac{1}{2} \mu \beta) } 
\ee
(compare the last line of \rf{proof-gamma}). 
If $\alpha, \beta \in {\cal R}_+$, then all four terms at the right-hand side are positive. Therefore, the left-hand side is positive which means that $|\mu \alpha\oplus\beta| < 2$. 
 \end{Proof}

The set ${\cal R}$ of all regressive functions 
is not closed with respect to $\oplus$. Indeed, suppose that a regressive $\alpha$ is given. The formula $\mu^2 \alpha \beta = - 4$ uniquely determines $\beta$ which also has to be regressive. However, in this case $\alpha\oplus\beta$ becomes infinite.

\ods
\begin{lem} \label{lem-albet}
\be
  x^\Delta (t) = \beta (t) x (t) \quad \Longleftrightarrow \quad 
 x^\Delta (t) = \alpha (t) \av{x (t)}  \ ,  
\ee
where $\beta$ is given by \rf{betal} and  
\be  \label{avg}
 \av{x (t)}  := \frac{ x (t) + x (\sigma (t)) }{2} \ . 
\ee
\end{lem}

\begin{Proof} Direct computation shows
\[
 x^\Delta  = \alpha  \av{x }  \quad \Longleftrightarrow \quad 
2 x^\Delta (1 + \frac{1}{2} \mu \beta) = \beta (x + x^\sigma) 
\quad \Longleftrightarrow \quad 
 x^\Delta  = \beta x \ ,
\]
where first we used  \rf{albet}, and then we substituted $x^\sigma = x + \mu x^\Delta$. 
\end{Proof}

\ods

\begin{Th}  \label{Th-Cau-Del}
The exponential function $E_\alpha (t, t_0)$, defined by \rf{E-for}, 
is a unique solution of the following Cauchy problem:
\be \label{xca}
       x^\Delta (t) = \alpha (t) \av{x (t)}  \ ,  \quad  x (t_0) = 1 \ , 
\ee
where $\alpha$ is regressive rd-continuous function and  $\av{x (t)}$ is defined \rf{avg}. 

\end{Th}

\begin{Proof} 
By Lemma~\ref{lem-albet} the initial value problem \rf{xca} is equivalent to 
\be \label{xcb}
   x^\Delta (t) = \beta (t) x (t) \ , \quad x (t_0) = 1 \ .
\ee
Taking into account that $\mu$ and $\alpha$ are rd-continuous, we have $\beta$ rd-continuous (compare Theorem~\ref{Ee}). Therefore, we can use Hilger's theorem concerning the problem \rf{xcb}, stating that its unique solution is $e_\beta (t, t_0)$.  
Finally, we use Theorem~\ref{Ee} once more. 
\end{Proof}

In the discrete case the equation \rf{xca} (treated as a numerical scheme) can be interpreted either as the trapezoidal rule, implicit midpoint rule or the discrete gradient method \cite{LaG,MQR2,HLW}. These implicit  methods are more accurate than the explicit (forward) Euler scheme (which is related to the equation $x^\Delta = \alpha x$) and can 
preserve more qualitative, geometrical and physical characteristics of the considered differential equations, e.g., integrals of motion \cite{HLW}.

\subsection{New definitions of hyperbolic and trigonometric functions}

Our definition of improved hyperbolic and trigonometric functions follows in a natural way (similarly as in the continuous case) from the definition of the exponential function.

\begin{Def}  \label{def-hyp}
Cayley-hyperbolic functions on a time scale are defined by
\be \label{def-chsh}
{\rm Cosh}_\alpha (t) := \frac{E_\alpha (t) + E_{-\alpha} (t) }{2} \ , \quad {\rm Sinh}_\alpha (t) := \frac{E_\alpha (t) - E_{-\alpha} (t) }{2} \ , 
\ee
where the exponential function $E_\alpha$ is defined by \rf{E-for}. 
\end{Def}

\begin{Def}  \label{def-trig}
Cayley-trigonometric functions on a time scale are defined by
\be  \label{cos-sin}
{\rm Cos}_\omega (t) := \frac{E_{i \omega} (t) + E_{- i \omega} (t) }{2} \ , \quad 
{\rm Sin}_\omega (t) := \frac{E_{i \omega} (t) - E_{- i \omega} (t) }{2 i} \ .
\ee
In other words, ${\rm Cos}_\omega (t) = {\rm Cosh}_{i \omega} (t)$, 
\  $i \, {\rm Sin}_\omega (t) = {\rm Sinh}_{i \omega} (t)$. 
\end{Def}

Properties of our hyperbolic and trigonometric functions are identical, or almost identical, as in the continuous case. Note that below we often use the notation \rf{avg}.

\begin{Th} Cayley-hyperbolic functions satisfy
\be \label{Pyth-hyp}
  {\rm Cosh}_\alpha^2 (t) - {\rm Sinh}_\alpha^2 (t) = 1 \ , 
\ee 
\be  \label{der-chsh}
{\rm Cosh}_\alpha^\Delta (t) = \alpha (t) \av{ {\rm Sinh}_\alpha (t)} \ , \qquad
{\rm Sinh}_\alpha^\Delta (t) = \alpha (t) \av{ {\rm Cosh}_\alpha (t)} \ .
\ee
\end{Th}

\begin{Proof} By Theorem~\ref{Th-exp-properties} we have $E_\alpha (t) E_{-\alpha} (t) = 1$. This is sufficient to directly verify the identity \rf{Pyth-hyp}. By Theorem~\ref{Th-Cau-Del} we have 
$E_\alpha^\Delta (t) = \alpha (t) \av{ E_\alpha (t) }$ and    
 differentiating \rf{def-chsh} we get \rf{der-chsh}.
\end{Proof}

\begin{Th} Cayley-trigonometric functions are real-valued for real $\omega$ and satisfy
\be \label{Pyth-trig}
{\rm Cos}_\omega^2 (t) + {\rm Sin}_\omega^2 (t) = 1 \ , 
\ee
\be  \label{der-cossin}
{\rm Cos}_\omega^\Delta (t) = - \omega (t) \av{ {\rm Sin}_\omega (t)} \ , \qquad
{\rm Sin}_\omega^\Delta (t) = \omega (t) \av{ {\rm Cos}_\omega (t)} \ .
\ee
\end{Th}

\begin{Proof} By Theorem~\ref{Th-exp-properties} we have $\overline{E_{i\omega } (t)} =  E_{-i \omega} (t)$. Thus the reality of trigonometric functions follows directly from Definition~\ref{def-trig}. Using $E_{i\omega } (t) E_{-i \omega} (t) = 1$ we get the Pythagorean identity \rf{Pyth-trig} starting from \rf{cos-sin}. Derivatives \rf{der-cossin} can be obtained by straightforward differentiation of the exponential functions using Theorem~\ref{Th-Cau-Del}.  
\end{Proof}

\begin{Th}  \label{Th-imag}
The function $\C \ni \alpha \rightarrow E_\alpha (t) $ maps the imaginary axis into 
the unit circle, i.e., \ ${\rm Re} \alpha (t) \equiv 0 \ \Longrightarrow \ |E_\alpha (t)| \equiv 1$.  
\end{Th}

\begin{Proof} We compute
\be
|E_\alpha (t)|^2 = E_\alpha (t) \overline{ E_\alpha (t)} = E_\alpha (t) E_{\bar \alpha} (t) = E_{\alpha\oplus{\bar \alpha}} \ ,
\ee
where we used twice Theorem~\ref{Th-exp-properties}. From the formula \rf{oplusmu} we see immediately that
\[
  \alpha \oplus \bar\alpha = 0 \quad \Longleftrightarrow \quad \bar\alpha = - \alpha \ ,
\]
i.e., $\alpha = i\omega$, where $\omega (t) \in \R$ for $t \in \T$. Therefore, $E_{\alpha\oplus\bar\alpha} (t) \equiv 1$ for imaginary $\alpha$. 
\end{Proof}

 Theorem~\ref{Th-imag} is crucial for a proper definition of trigonometric functions. Indeed, $|E_{i\omega}| = 1$  is equivalent to 
$ ({\rm Re} E_{i\omega})^2 + ({\rm Im} E_{i\omega})^2 = 1$ and we can identify $ {\rm Re} E_{i\omega}$ and ${\rm Im} E_{i\omega}$ with trigonometric functions.

At the end of this section we present second order dynamic equations satisfied by hyperbolic and trigonometric functions with constant $\alpha$ and $\omega$, respectively. 

\begin{lem}  \label{lem-av-com}
Averaging commutes with delta differentiation, i.e., 
\be  \label{av-del-com}
\av{x (t)}^\Delta = \av{ x^\Delta (t) } \ .
\ee
\end{lem}
\begin{Proof} \quad $ \dis
\av{x (t)}^\Delta = \frac{1}{2} \left( x (t) + x(t^\sigma) \right)^\Delta = \frac{1}{2} \left( x^\Delta (t) + x^\Delta (t^\sigma) \right) = \av{ x^\Delta (t) } \ .
$
\end{Proof}

\begin{prop} \label{prop-hyp-osc}
If \ $\alpha (t) = \const$, then improved hyperbolic functions on a time scale, ${\rm Cosh}_\alpha$ and ${\rm Sinh}_\alpha$,  satisfy the equation
\be 
 x^{\Delta\Delta} (t) = \alpha^2 \av{\av{x (t)}} \ , 
\ee
where \ $\av{\av{ x (t) }} \equiv \frac{1}{4}  \left( x + x^\sigma + x^{\sigma\sigma} \right)$ \ and \ $x^\sigma (t) :=  x (t^\sigma)$. 
\end{prop}

\begin{Proof} follows immediately from \rf{der-chsh} and \rf{av-del-com}. Indeed,  
\[
 {\rm Cosh}_\alpha^{\Delta\Delta} (t) = \alpha \av{\rm Sinh_\alpha (t) }^\Delta = \alpha \av{ \alpha \av{{\rm Cosh_\alpha (t) } } } = \alpha^2 \av{\av{ {\rm Cosh_\alpha (t) }} } \ .
\]
The same calculation can be done for $\rm Sinh_\alpha^{\Delta\Delta}$.  
\end{Proof}

\begin{prop} \label{prop-osc}
If \ $\omega (t) = \const$, then improved trigonometric functions on a time scale, ${\rm Cos}_\omega$ and ${\rm Sin}_\omega$, satisfy the equation
\be  \label{oscyl}
 x^{\Delta\Delta}  + \omega^2 \av{\av{x (t)}} = 0 \ .
\ee
\end{prop}

\begin{Proof} by straightforward calculation, compare the proof of Proposition~\ref{prop-hyp-osc}.  \end{Proof}

\section{Exact special functions on time scales}
\label{sec-exact}

The simplest (almost trivial) way to construct time scales analogues of special functions is to take their {\it exact} values. 

\begin{Def}
Given a function $f : \R \rightarrow \C$ we define its {\em exact} analogue  $\tilde f : \T \rightarrow \C$ as $\tilde f := f|_{\T}$, i.e., 
\be
    \tilde f (t) := f (t) \quad ({\rm for}\ t \in \T) \ .
\ee
\end{Def}

Although the path $f \rightarrow \tilde f$ is obvious and unique, the inverse correspondence may cause serious problems.   Usually, it is not  easy to find or to indicate the most appropriate (or most natural) real function corresponding to a given function on $\T$.

\subsection{Exact exponential function on time scales}

The continuous exponential function is given by $e_a (t) = \exp \int_{t_0}^t a (\tau) d \tau$, where $a : \R \rightarrow \C$ is a given function. 
A non-trivial question is to choose a function $a = a(t)$, provided that  we intend to define an exact exponential corresponding to a given  function   $\alpha : \T \rightarrow \C$.  In general, the choice seems to be highly non-unique.

In this paper we confine ourselves to the simplest case, $\alpha = \const$, when it is natural to take \ $a (t) = \alpha = \const$.  

\begin{Def}
The exact exponential function $E^{ex}_\alpha (t, t_0)$ (where $\alpha$ is a complex constant and $t, t_0 \in \T$) is defined by 
$E^{ex}_\alpha (t, t_0) := e^{\alpha (t-t_0)}$. 
\end{Def}

\begin{Th}
The exact exponential function $E^{ex}_\alpha (t, t_0)$ 
satisfies the dynamic equation
\be  \label{dyn-ex}
  x^\Delta (t)  = \alpha \psi_{\alpha} (t)    \av{x (t)} \ , 
\qquad x (t_0) = 1 \ ,
\ee
where $\psi_\alpha (t) = 1$ for right-dense points and  
\be  \label{psi}
  \psi_\alpha (t) = \frac{2}{\alpha \mu (t) }{\tanh \frac{\alpha \mu (t)}{2} }
\ee
for right-scattered points. 
\end{Th}

\begin{Proof} For right-dense $t$ the equation \rf{dyn-ex} reduces to $x^\Delta = \alpha x$. Then, still assuming $t$ right-dense, we compute 
\[
(e^{\alpha (t-t_0)})^\Delta = \frac{d}{dt} e^{\alpha (t-t_0)} = \alpha e^{\alpha (t-t_0)} \ ,
\]
i.e., $e^{\alpha (t-t_0)}$ satisfies the equation $x^\Delta = \alpha x$. 
For right-scattered $t$ we have: 
\[  \ba{l} \dis
  E^{ex}_\alpha (t^\sigma, t_0) - E^{ex}_\alpha (t, t_0) = e^{\alpha (t^\sigma - t_0) } - e^{\alpha (t - t_0)} =  e^{\alpha (t - t_0)} \left( e^{\alpha \mu} - 1 \right) \ , \\[2ex]
E^{ex}_\alpha (t^\sigma,t_0) + E^{ex}_\alpha (t, t_0) = e^{\alpha (t^\sigma - t_0) } + e^{\alpha (t - t_0)}  =  e^{\alpha (t - t_0)} \left( e^{\alpha \mu} + 1 \right)  \ ,
\ea \]
and, therefore,
\be
 E^{ex}_\alpha (t^\sigma, t_0) - E^{ex}_\alpha (t, t_0) = \tanh\frac{\alpha\mu}{2} \left( E^{ex}_\alpha (t^\sigma, t_0) + E^{ex}_\alpha (t, t_0)  \right) \ ,
\ee
which is equivalent to \rf{dyn-ex}. The initial condition is obviously satisfied. 
\end{Proof}

The equation \rf{dyn-ex} can be interpreted as the dynamic equation \rf{xca} with a modified $\alpha$ (i.e., $\alpha \rightarrow \alpha \psi_\alpha $). Another interpretation can be obtained by a modification of the delta derivative. We define
\be
  x^{\Delta'_\alpha} (t) := \lim_{\stackrel{\dis s \rightarrow t}{s \neq \sigma (t)} } 
\frac{ x (t^\sigma) - x (s)}{ \delta_\alpha (t^\sigma - s)  }
\ee
where $\delta_\alpha$ is a function  given by 
\be
\delta_\alpha (\mu) : = \frac{2}{\alpha} \tanh \frac{ \alpha \mu}{2} \ .
\ee

\begin{lem}
\be  \label{delt1}
   x^\Delta (t) = \psi_\alpha (t) x^{\Delta'_\alpha} (t)
\ee
\end{lem}

\begin{Proof} We compute
\[  \lim_{\stackrel{\dis s \rightarrow t}{s \neq \sigma (t)} }\frac{ x (t^\sigma) - x (s)}{ t^\sigma - s } = 
\lim_{\stackrel{\dis s \rightarrow t}{s \neq \sigma (t)} } 
\frac{ x (t^\sigma) - x (s)}{ \delta_\alpha (t^\sigma - s)  }   \lim_{\stackrel{\dis s \rightarrow t}{s \neq \sigma (t)} } 
\frac{ \delta_\alpha (t^\sigma - s)  }{ t^\sigma - s} = x^{\Delta'_\alpha} (t)   \psi_\alpha (t)  \ ,
\]
which yields \rf{delt1}. 
\end{Proof}

\begin{cor} The equation \rf{dyn-ex}, satisfied by the exact exponential function, can be rewritten as  
\be
x^{\Delta'_\alpha} (t) = \alpha \av{ x (t) } \ , \qquad x (t_0) = 1 \ .
\ee
\end{cor}

The equation \rf{dyn-ex} is the exact discretization of the equation $\dot x = \alpha  x$. In general, 
by the exact discretization of an ordinary differential equation $\dot x = f (x)$, where $x (t) \in \R^N$, we mean the difference equation $X_{n+1} = F (X_n)$, where $X_n \in \R^N$, such that $X_n = x (t_n)$. Any equation has an implicit exact discretization (provided that the solution exists). It is worthwhile to point out that all linear ordinary differential equations with constant coefficients have explicit exact discretizations \cite{Po,Mic,Ag}.

\subsection{Exact hyperbolic and trigonometric functions on time scales} 

In order to simplify notation we confine ourselves to $t_0 = 0$. All results can be obviously extended on the general case.

\begin{Def} Given a real constant $\alpha$, exact hyperbolic  functions on a time scale $\T$ are defined by
\be  \ba{l} \dis
\cosh^{ex}_\alpha (t) =  \frac{ E^{ex}_{\alpha} (t) + E^{ex}_{-\alpha} (t)}{2} = \cosh \alpha  t \ , \\[2ex] \dis
\sinh^{ex}_\alpha (t) =  \frac{ E^{ex}_{\alpha} (t) - E^{ex}_{-\alpha} (t)}{2} = \sinh\alpha t \ .
\ea \ee
\end{Def}

 We point out that in this definition $t\in\T$. The same remark concerns the next definition. 

\begin{Def} Given a real constant $\omega$, exact trigonometric functions on a time scale $\T$ are defined by
\be  \ba{l} \dis
\cos^{ex}_\omega (t) =  \frac{ E^{ex}_{i\omega} (t) + E^{ex}_{-i\omega} (t)}{2} = \cos \omega t \ , \\[2ex] \dis
\sin^{ex}_\omega (t) =  \frac{ E^{ex}_{i\omega} (t) - E^{ex}_{-i\omega} (t)}{2 i} = \sin\omega t \ .
\ea \ee
\end{Def}

It turns out that dynamic equations satisfied by exact trigonometric and hyperbolic functions are rather awkward in the case of arbitrary time scales. These equations simplify considerably under asumption of constant graininess. Therefore, from now on, we confine ourselves to the case $\mu (t) = \const$. 

Moreover, we observe that the function $\psi_\alpha$ (defined by \rf{psi}) is symmetric. In particular, $\psi_{i\omega} (t) = \psi_{-i\omega} (t)$. We denote: 
\be  \label{phi}
\phi ( x ) := \frac{2}{ x } \tan \frac{ x}{2} \quad ({\rm for} \ x \neq 0) \ , \quad \phi (0) = 1\ .
\ee
Therefore, for $\mu (t) = \const$ we have
\be
 \psi_{i \omega} (t) = \phi (\omega\mu) = \const \ .
\ee

\begin{Th}
We assume that the time scale $\T$ has constant graininess $\mu$. 
Then exact trigonometric functions on $\T$, i.e., $\cos^{ex}_\omega$ and $\sin^{ex}_\omega$, 
satisfy the dynamic equation
\be  \label{osc-avav}
 x^{\Delta\Delta} (t) + \omega^2 \phi^2 (\omega\mu) \, \av{\av{ x (t) }} = 0 \ ,
\ee
which is equivalent to 
\be  \label{osc-sinc-ex}
 x^{\Delta \Delta} (t) + \omega^2 \left({\rm sinc} \frac{\omega\mu}{2} \right)^2  x (t^\sigma ) = 0 \ ,
\ee
where
\be  \label{sinc}
{\rm sinc} (x) := \frac{\sin x}{x}  \qquad ({\rm for} \ \ x \neq 0) \ , \qquad {\rm sinc}( 0 ) := 1 \ . 
\ee
\end{Th}

\begin{Proof} By Theorem~\ref{Th-Cau-Del} we have
\be
  ( E^{ex}_{\pm i\omega} )^\Delta = \pm i\omega \psi_{\pm i\omega} (t) \av{ E^{ex}_{\pm i\omega} } \ .
\ee
Therefore, taking also \rf{phi} into account, we get
\be  \ba{l} \label{delta-sca}
(\cos^{ex}_\omega (t))^\Delta = \frac{1}{2} \left( E^{ex}_{i\omega} (t) + E^{ex}_{-i\omega} (t) \right)^\Delta = - \omega \phi (\omega\mu) \av{ \sin^{ex}_\omega (t) } \ , 
\\[3ex]
(\sin^{ex}_\omega (t))^\Delta = \frac{1}{2i} \left( E^{ex}_{i\omega} (t) - E^{ex}_{-i\omega} (t) \right)^\Delta =  \omega \phi (\omega\mu) \av{ \cos^{ex}_\omega (t) } \ .
\ea \ee  
We notice that the second formula is equivalent to \rf{delta-sin-ex}.  Applying Lemma~\ref{lem-av-com} to \rf{delta-sca}, we obtain 
\be \ba{l}
( \cos^{ex}_\omega (t) )^{\Delta \Delta} = - \omega^2 \phi^2 (\omega\mu) \av{\av{ \cos^{ex}_\omega (t) }} \ ,  \\[2ex]
( \sin^{ex}_\omega (t) )^{\Delta \Delta} = - \omega^2 \phi^2 (\omega\mu) \av{\av{ \sin^{ex}_\omega (t) }} \ ,
\ea \ee
i.e., we have \rf{osc-avav}. The equivalence between \rf{osc-avav} and \rf{osc-sinc-ex} is obvious for $\mu = 0$. In the case $\mu \neq 0$ we use trigonometric identities:
\be \ba{l} \dis
\sin(\omega t) + \sin(\omega t + 2 \omega \mu) + 2 \sin(\omega t + \omega\mu) = 2 (1 + \cos\omega\mu) \sin(\omega t + \omega \mu) , \\[2ex] \dis 
\cos(\omega t) + \cos(\omega t + 2 \omega \mu) + 2 \cos(\omega t + \omega\mu) = 2 (1 + \cos\omega\mu) \cos(\omega t + \omega \mu) , 
\ea \ee
obtaining
\be \ba{l} \dis
\av{\av{ \sin^{ex}_\omega (t) }} = \left( \cos\frac{\omega\mu}{2} \right)^2 \sin(\omega t + \omega\mu) \ ,  \\[2ex] \dis
\av{\av{ \cos^{ex}_\omega (t) }} = \left( \cos\frac{\omega\mu}{2} \right)^2 \cos(\omega t + \omega\mu) \ .
\ea \ee
Taking into account that \ ${\rm sinc} \frac{\omega\mu}{2} = \phi (\omega\mu)  \cos \frac{\omega\mu}{2} $ \  and \ $t^\sigma = t + \mu$, we get
\be
  \phi^2 (\omega\mu) \av{\av{ x (t) }} = \left( {\rm sinc} \frac{\omega\mu}{2} \right)^2 x (t^\sigma) \ ,
\ee
where $x (t)$ is an arbitrary linear combination of $\sin^{ex}_{\omega}$ and $\cos^{ex}_{\omega}$, which ends the proof. 
\end{Proof}

The exact discretization of the harmonic oscillator equation  leads to  
another modification of the delta derivative (see 
\cite{Ci-oscyl,CR-ade}), 
\be  \label{deltabis}
  x^{\Delta''_\omega} (t) = \lim_{\stackrel{\dis s \rightarrow t}{s \neq \sigma (t)} } 
\frac{ x (t^\sigma) - x (s) \cos  (\omega t^\sigma - \omega s  )}{\omega^{-1} \sin (\omega t^\sigma - \omega s)} \ .
\ee

In order to avoid infinite values of $x^{\Delta''_\omega}$ we assume $|\omega \mu (t) | < \pi$. All positively regressive constant functions $\omega$ (see Definition~\ref{Def-posit}) obviously satisfy this requirement.

\begin{prop}
If $x = x (t)$ satisfies $\ddot x + \omega^2 x = 0$ for $t\in \R$, then 
\be   \label{pochodne} 
  ( x (t) |_{t\in \T} )^{\Delta''_\omega } = \dot x (t) |_{t\in \T} \ .
\ee 
\end{prop}

\begin{Proof}
By assumption,  $x (t) = A \cos\omega t + B \sin \omega t$. Then 
\[
  x (t^\sigma) = \cos \omega\mu \ (A \cos\omega t + B \sin \omega t) + \sin\omega\mu \left(  B \cos\omega t - A \sin\omega t \right) \ ,
\]
because $t^\sigma = t+\mu$. By direct computation we verify 
\[
 x (t^\sigma ) - x (s) \cos (\omega t^\sigma - \omega s) = ( B \cos \omega s - A \sin \omega s ) \sin (\omega t^\sigma - \omega s) \ . 
\]
Therefore,
\[
 x^{\Delta''_\omega} (t) = \lim_{\stackrel{\dis s \rightarrow t}{s \neq \sigma (t)} } 
\frac{ x (t^\sigma) - x (s) \cos  (\omega t^\sigma - \omega s  )}{\omega^{-1} \sin (\omega t^\sigma - \omega s)}  = \omega (B \cos \omega t - A \sin \omega t ) = \dot x (t) ,
\]
which ends the proof. 
\end{Proof}

\begin{lem} \label{lem-delbis}
\be  \label{delbis}
x^\Delta (t) = {\rm sinc} (\omega\mu) \ x^{\Delta''_\omega} (t) - \frac{1}{2} \mu \omega^2   \left( {\rm sinc} \frac{\omega \mu }{2} \right)^2 x (t) \ .
\ee
\end{lem}

\begin{Proof} Substituting the definition \rf{deltabis} into the right-hand side of the formula \rf{delbis}, and taking into account that 
 \be \ba{l} \dis
\lim_{\stackrel{\dis s \rightarrow t}{s \neq \sigma (t)} } \frac{ \sin (\omega t^\sigma - \omega s)  }{\omega t^\sigma - \omega s } = {\rm sinc (\omega\mu) }  \ , \\[3ex] \dis
  \lim_{\stackrel{\dis s \rightarrow t}{s \neq \sigma (t)} } \frac{ x (s) ( 1 - \cos (\omega t^\sigma - \omega s) )}{ t^\sigma - s} = \frac{1}{2} \mu \omega^2   \left( {\rm sinc} \frac{\omega \mu }{2} \right)^2 x (t)\ ,
\ea \ee
we obtain the left-hand side of  \rf{delbis}, compare \rf{delta}. 
\end{Proof}

\begin{prop} The equation \rf{osc-sinc-ex}, satisfied by exact trigonometric functions, can be rewritten as  
\be  \label{ex-osc}
  x^{ \Delta''_\omega \Delta''_\omega} (t) + \omega^2 x (t) = 0 \ . 
\ee
\end{prop}

\begin{Proof}
The formula  \rf{delbis} can be rewritten as
\be
 x^\Delta (t) = {\rm sinc} \frac{\omega\mu}{2} \left( \cos \frac{\omega\mu}{2} \ x^{\Delta''_\omega} (t)  - \omega \sin \frac{\omega\mu}{2} \ x (t) \right) \ ,
\ee
and, therefore,
\[
 x^{\Delta \Delta} (t) = {\rm sinc}^2 \frac{\omega\mu}{2} 
\left(  \cos^2 \frac{\omega\mu}{2} \ x^{\Delta''_\omega \Delta''_\omega} (t)  - \omega (\sin\omega\mu) \, x^{\Delta''_\omega} (t) + \omega^2 \sin^2 \frac{\omega\mu}{2} \, x (t) \right) .
\]
Moreover, $x (t^\sigma) = x (t) + \mu x^\Delta (t)$. Hence
\be
x (t^\sigma) = (\cos\omega\mu) \ x (t) + \frac{\sin \omega\mu}{\omega} \ x^{\Delta''_\omega} (t) \ . 
\ee
Thus we can easily verify that
\be
x^{\Delta\Delta} (t) + \omega^2 {\rm sinc}^2 \frac{\omega\mu}{2} x (t) \equiv 
{\rm sinc}^2 \frac{\omega\mu}{2} \cos^2 \frac{\omega\mu}{2} \left( 
x^{ \Delta''_\omega \Delta''_\omega} (t) + \omega^2 x (t) \right)  , 
\ee
which ends the proof. 
\end{Proof}

The exact discretization of the harmonic oscillator equation $\ddot x + \omega^2 x = 0$ was discussed in detail in \cite{Ci-oscyl}. In particular, we presented there the discrete version of the formula \rf{ex-osc} and related results.

\section{Conclusions and future directions}

We proposed two different new approaches to the construction of exponential, hyperbolic and trigonometric functions on time scales. 

The resulting functions preserve most of the qualitative properties of the corresponding continuous functions. In particular, Pythagorean trigonometric identities hold exactly on any time scale. Dynamic equations satisfied by Cayley-motivated functions have a natural similarity to the corresponding diferential equations. 

The first approach is based on the Cayley transformation. 
It has important advantages because simulates better the behaviour of the exponential function, both qualitatively (e.g., the Cayley-exponential function maps the imaginary axis into the unit circle) and quantitatively, because 
\be
  \frac{1 + \frac{1}{2} \mu \alpha}{1 - \frac{1}{2} \mu \alpha} = 
1 + \alpha \mu + \frac{1}{2} (\alpha \mu)^2 + \ldots
\ee
i.e., this factor approximates $\exp (\alpha \mu)$ up to second order terms, while $1 + \mu \alpha$ or $(1 - \mu \alpha)^{-1}$ are only first order approximations.

Our approach has some disadvantages, as well. The promising notion of complete delta differentiability \cite{BG-part,Ci-pst} becomes very difficult or impossible to apply (because integral curves of our new dynamic systems on time scales, like \rf{xca}, are not completely delta differentiable). Moreover, 
dynamic equation on time scales become implicit (and equations of standard delta calculus on time scales are explicit). However, nabla calculus is also implicit, and it is very well known that implicit numerical finite difference schemes have better properties than explicit schemes, see for instance \cite{HLW}. 

The second approach consists in exact discretization. 
In this paper we confined ourselves to the simplest case, i.e., to the exponential function $E_\alpha (t)$ with $\alpha = \const$. It would be interesting to define and study exact exponential functions for a larger class of functions $\alpha$. We leave it as an interesting open problem. 
Other problems are associated with finding dynamic systems which correspond to exact discretizations. 
In the case of linear equations this subject is well known, but more general results are diffcult to be obtained. 

We point out that definitions of elementary functions on time scales are  not unique. We presented some arguments in advantage of our definitions but in principle one can develop several different theories of elementary and special functions on time scales. It seems important to develop and understand different approaches, and, if possible, to find a ``vocabulary'' to translate results. 
Similarly, one can develop several different theories of dynamic systems on time scales closely related to different numerical finite difference schemes.
 For instance, the standard delta calculus corresponds to forward (explicit) Euler scheme, the nable calculus corresponds to the implicit Euler scheme, and my proposition is related to the trapezoidal rule (and to the discrete gradient methods). Therefore there are no unique ``natural'' time scales analogues of dynamic systems. One can choose among many possibilities, including the above three approaches and the exact discretization (which explicitly exists only for a very limited class of differential equations). Another  promising possibility is the so called ``locally exact'' discretization \cite{Ci-oscyl,CR-grad}.  

The definitions presented in our paper seem to be entirely new as far as time scales are concerned but their discrete analogues ($\T = \ep \Z$) have been used since a long time. After completing this work I found a lot of references where rational or Cayley-like forms of the exponential function are  used in the discrete case, see for instance \cite{Is-Cay,Fer,Duf,ZD,DJM1,NQC,Mer,BMS}. It would be interesting to specify our results to the quantum calculus case \cite{KC}, where the approach presented in this paper probably has not been applied yet. 

It is obvious that the proposed modification of the basic definitions should have an essential influence on many branches of the time scales calculus, including the theory of dynamic equations \cite{ABOP,BP-I}, Hamiltonian systems \cite{ABR}, and the Fourier and Laplace transforms \cite{BG-Lap}. 

Finally, we notice that exponential functions on time scales are defined (here and in other papers) on real time scales. It would be important to extend these definitions on the complex domain which is so natural for continuous exponential functions. Such extensions are well known in the discrete case and in the quantum calculus, see for instance \cite{Fer,Duf,Mer,BMS}. 
\ods

{\it Acknowledgements.} I am grateful to Stefan Hilger for encouragement and sending me the paper \cite{Hi-spec}, to Maciej Nieszporski for turning my attention on Mercat's and Nijhoff's papers \cite{NQC,Mer}, and to Adam Doliwa for the reference \cite{DJM1}. Discussions with Maciej Nieszporski concerning discretizations of Lax pairs, see \cite{CMNN}, turned out to be useful also in the context of this paper.


\end{document}